\documentclass{article}

\usepackage{amsmath,amssymb,bm}
\usepackage{color,url,cite}
\usepackage{arydshln}

\usepackage[left=2.0cm,right=2.0cm,top=2.0cm,bottom=2.0cm]{geometry}

\newtheorem{theorem}{Theorem}
\newtheorem{lemma}{Lemma}
\newtheorem{corollary}{Corollary}
\newenvironment{proof}[1][Proof]{\par\noindent\textit{#1}: }{\hfill$\blacksquare$\vskip 0.5\baselineskip}

\begin{document}

\title{Evaluating Two Determinants}
\author{Shujun Li\\\url{http://www.hooklee.com}}
\date{\today}

\maketitle

\begin{abstract}
This article evaluates the determinants of two classes of special
matrices, which are both from a number theory problem. Applications
of the evaluated determinants can be found in
[arXiv:math.NT/0509523, 2005].

Note that the two determinants are actually special cases of
Theorems 20 and 23 in [arXiv:math.CO/9902004], respectively. Since
this paper does not provide any new results, it will not be
published anywhere.
\end{abstract}

\section{The Determinants}

\begin{theorem}\label{theorem:det-powers}
Assume $m\geq 1$. Given a $2m\times 2m$ matrix
$\bm{A}=\left[\begin{matrix}\bm{A}_1\\\bm{A}_2\end{matrix}\right]$,
where $\bm{A}_1=[X_i^{j-1}]_{1\leq i\leq m \atop 1\leq j\leq 2m}$
and $\bm{A}_2=[jX_i^{j-1}]_{1\leq i\leq m \atop 1\leq j\leq 2m}$,
i.e.,
\[
\bm{A}=\left[\begin{array}{cccc:ccc}%
1 & X_1 & \cdots & X_1^{m-1} & X_1^m & \cdots & X_1^{2m-1}\\
1 & X_2 & \cdots & X_2^{m-1} & X_2^m & \cdots & X_2^{2m-1}\\
\vdots & \vdots & \ddots & \vdots & \vdots & \ddots & \vdots\\
1 & X_m & \cdots & X_m^{m-1} & X_m^m & \cdots & X_m^{2m-1}\\
\hdashline
1 & 2X_1 & \cdots & mX_1^{m-1} & (m+1)X_1^m & \cdots & 2mX_1^{2m-1}\\
1 & 2X_2 & \cdots & mX_2^{m-1} & (m+1)X_2^m &
\cdots & 2mX_2^{2m-1}\\
\vdots & \vdots & \ddots & \vdots & \vdots & \ddots & \vdots\\
1 & 2X_m & \cdots & mX_m^{m-1} & (m+1)X_m^m & \cdots & 2mX_m^{2m-1}
\end{array}\right].
\]
Then, $|\bm{A}|=(-1)^{\frac{m(m-1)}{2}}\prod_{i=1}^mX_i\prod_{1\leq
i<j\leq m}(X_j-X_i)^4$.
\end{theorem}

\begin{corollary}\label{corollary:det-powers}
Assume $m\geq 1$. Given a $2m\times 2m$ matrix
$\bm{A}=\left[\begin{matrix}\bm{A}_1\\\bm{A}_2\end{matrix}\right]$,
where $\bm{A}_1=[X_i^{j+1}]_{1\leq i\leq m \atop 1\leq j\leq 2m}$
and $\bm{A}_2=[(j+1)X_i^j]_{1\leq i\leq m \atop 1\leq j\leq 2m}$.
Then,
$|\bm{A}|=(-1)^{\frac{m(m-1)}{2}}\prod_{i=1}^mX_i^4\prod_{1\leq
i<j\leq m}(X_j-X_i)^4$.
\end{corollary}
\begin{proof}
Factoring out $X_i^2$ from each row of $\bm{A}_1$ and factoring out
$X_i$ from each row of $\bm{A}_2$, one has
$\bm{A}_1^{(1)}=[X_i^{j-1}]_{1\leq i\leq m \atop 1\leq j\leq 2m}$
and $\bm{A}_2^{(1)}=[(j+1)X_i^{j-1}]_{1\leq i\leq m \atop 1\leq
j\leq 2m}$. Then, for $i=1\sim m$, subtracting row $i$ of
$\bm{A}_1^{(1)}$ from row $i$ of $\bm{A}_2^{(1)}$, one has
$\bm{A}_2^{(2)}=[jX_i^{j-1}]_{1\leq i\leq m \atop 1\leq j\leq 2m}$.
From Theorem \ref{theorem:det-powers}, one immediately gets
\[
|\bm{A}|=\prod_{i=1}^mX_i^3\left((-1)^{\frac{m(m-1)}{2}}\prod_{i=1}^mX_i\prod_{1\leq
i<j\leq
m}(X_j-X_i)^4\right)=(-1)^{\frac{m(m-1)}{2}}\prod_{i=1}^mX_i^4\prod_{1\leq
i<j\leq m}(X_j-X_i)^4.
\]
\end{proof}

\begin{theorem}\label{theorem:det-binom-powers}
Assume $m\geq 1,n\geq l\geq 1$ and $\bm{A}$ is a block-wise
$ml\times ml$ matrix as follows:
\[
\bm{A}=\left[\begin{matrix}%
\bm{A}_1\\
\bm{A}_2\\
\vdots\\
\bm{A}_m
\end{matrix}\right],
\]
where for $i=1\sim m$,
\[
\bm{A}_i=\left[\binom{n+j-1}{k-1}X_i^{j-1}\right]_{1\leq
j\leq ml \atop 1\leq k\leq l}=\left[\begin{matrix}%
\binom{n}{0} & \binom{n+1}{0}X_i & \cdots &
\binom{n+(ml-1)}{0}X_i^{ml-1}\\
\binom{n}{1} & \binom{n+1}{1}X_i & \cdots &
\binom{n+(ml-1)}{1}X_i^{ml-1}\\
\vdots & \vdots & \ddots & \vdots\\
\binom{n}{l-1} & \binom{n+1}{l-1}X_i & \cdots &
\binom{n+(ml-1)}{l-1}X_i^{ml-1}
\end{matrix}\right]_{ml\times l}.
\]
Then, $|\bm{A}|=\prod_{i=1}^mX_i^{\frac{l(l-1)}{2}}\prod_{1\leq
i<j\leq m}(X_j-X_i)^{l^2}$.
\end{theorem}

\section{The Proofs}

\subsection{Proof of Theorem \ref{theorem:det-powers}}

\begin{proof}
We use mathematical induction on $m$ to prove this theorem.

1) When $m=1$, $\bm{A}=\left[\begin{matrix}1 & X_1\\1 &
2X_1\end{matrix}\right]$. Directly calculating the determinant,
$|\bm{A}|=X_1=(-1)^{\frac{1(1-1)}{2}}X_1$.

2) Assume this theorem is true for $m-1\geq 1$, let us prove the
case of $m\geq 2$.

For $j=2\sim 2m$, subtracting $X_1$ times of column $(j-1)$ from
column $j$ of $\bm{A}$, one gets
\[
\left[\begin{array}{c:c}%
\begin{matrix}1\end{matrix} & \begin{matrix}0 & 0 & \cdots & 0\end{matrix}\\
\hdashline
\begin{matrix}1\\\vdots\\1\end{matrix} & \bm{A}_1^{(1)}\\
\hdashline
\begin{matrix}1\end{matrix}
& \begin{matrix}X_1 & X_1^2 & \cdots & X_1^{2m}\end{matrix}\\
\hdashline
\begin{matrix}1\\\vdots\\1\end{matrix} & \bm{A}_2^{(1)}
\end{array}\right],
\]
where $\bm{A}_1^{(1)}$ and $\bm{A}_2^{(1)}$ are both
$(m-1)\times(2m-1)$ matrices:
$\bm{A}_1^{(1)}=[X_{i+1}^{j-1}(X_{i+1}-X_1)]_{1\leq i\leq m-1 \atop
1\leq j\leq 2m-1}$ and
$\bm{A}_2^{(1)}=[jX_{i+1}^{j-1}(X_{i+1}-X_1)+X_{i+1}^j]_{1\leq i\leq
m-1 \atop 1\leq j\leq 2m-1}$. Apparently, $|\bm{A}|$ is equal to the
determinant of the following $(2m-1)\times(2m-1)$ matrix:
\[
\bm{A}^{(1)}=\left[\begin{matrix}%
\bm{A}_1^{(1)}\\
\begin{matrix}X_1 & X_1^2 & \cdots & X_1^{2m}\end{matrix}\\
\bm{A}_2^{(1)}
\end{matrix}\right].
\]
Moving the row matrix between $\bm{A}_1^{(1)}$ and $\bm{A}_2^{(1)}$
to the top of $\bm{A}^{(1)}$, one get another matrix $\bm{A}^{(2)}$
and has $|\bm{A}|=(-1)^{m-1}|\bm{A}^{(2)}|$. Factoring
$(X_{i+1}-X_1)$ out from each row of $\bm{A}_1^{(1)}$, one gets a
new sub-matrix $\bm{A}_1^{(3)}=[X_{i+1}^{j-1}]_{1\leq i\leq m-1
\atop 1\leq j\leq 2m-1}$ and
\[
\bm{A}^{(3)}=\left[\begin{matrix}%
\begin{matrix}X_1 & X_1^2 & \cdots & X_1^{2m}\end{matrix}\\
\bm{A}_1^{(3)}\\
\bm{A}_2^{(1)}
\end{matrix}\right].
\]
Apparently, $|\bm{A}|=(-1)^{m-1}\prod_{1\leq i\leq
m-1}(X_{i+1}-X_1)|\bm{A}^{(3)}|=(-1)^{m-1}\prod_{2\leq i\leq
m}(X_i-X_1)|\bm{A}^{(3)}|$. Then, for $i=1\sim m-1$, subtract row
$i$ of $\bm{A}_1^{(3)}$ multiplied by $X_{i+1}$ from row $i$ of
$\bm{A}_2^{(1)}$, one has a new sub-matrix
$\bm{A}_2^{(4)}=[jX_{i+1}^{j-1}(X_{i+1}-X_1)]_{1\leq i\leq m-1 \atop
1\leq j\leq 2m-1}$. Then, for $i=1\sim m-1$, factor out
$(X_{i+1}-X_1)$ from each row of $\bm{A}_2^{(4)}$, one has
$\bm{A}_2^{(5)}=[jX_{i+1}^{j-1}]_{1\leq i\leq m-1 \atop 1\leq j\leq
2m-1}$ and
\[
\bm{A}^{(5)}=\left[\begin{matrix}%
\begin{matrix}X_1 & X_1^2 & \cdots & X_1^{2m}\end{matrix}\\
\bm{A}_1^{(3)}\\
\bm{A}_2^{(5)}
\end{matrix}\right].
\]
Now, $|\bm{A}|=(-1)^{m-1}\prod_{2\leq i\leq
m}(X_i-X_1)^2|\bm{A}^{(5)}|$. Next, for $j=2\sim 2m-1$, subtract
$X_1$ times of the column $(j-1)$ from the column $j$ of
$\bm{A}^{(5)}$, one has
\[
\bm{A}^{(6)}=\left[\begin{matrix}%
\begin{matrix}X_1\end{matrix} & \begin{matrix}0 & \cdots & 0\end{matrix}\\
\begin{matrix}1\end{matrix} & \bm{A}_1^{(6)}\\
\begin{matrix}1\end{matrix} & \bm{A}_2^{(6)}
\end{matrix}\right],
\]
where $\bm{A}_1^{(6)}$ and $\bm{A}_1^{(6)}$ are both $(m-1)\times
(2m-2)$ sub-matrices:
$\bm{A}_1^{(6)}=[X_{i+1}^{j-2}(X_{i+1}-X_1)]_{1\leq i\leq m-1 \atop
1\leq j\leq 2m-2}$ and
$\bm{A}_2^{(6)}=[jX_{i+1}^{j-1}(X_{i+1}-X_1)+X_{i+1}^j]_{1\leq i\leq
m-1 \atop 1\leq j\leq 2m-2}$. Assuming
$\bm{A}^{(7)}=\left[\begin{matrix}\bm{A}_1^{(6)}\\\bm{A}_2^{(6)}\end{matrix}\right]$,
one has $|\bm{A}|=(-1)^{m-1}X_1\prod_{2\leq i\leq
m}(X_i-X_1)^2|\bm{A}^{(7)}|$. Then, $\bm{A}_1^{(6)}$ and
$\bm{A}_2^{(6)}$ can be processed in the same way as
$\bm{A}_1^{(1)}$ and $\bm{A}_2^{(1)}$, one can get
$\bm{A}^{(8)}=\left[\begin{matrix}\bm{A}_1^{(8)}\\\bm{A}_2^{(8)}\end{matrix}\right]$,
where $\bm{A}_1^{(8)}=[X_{i+1}^{j-1}]_{1\leq i\leq m-1 \atop 1\leq
j\leq 2m-2}$ and $\bm{A}_2^{(8)}=[jX_{i+1}^{j-1}]_{1\leq i\leq m-1
\atop 1\leq j\leq 2m-2}$. Now, $|\bm{A}|=(-1)^{m-1}X_1\prod_{2\leq
i\leq m}(X_i-X_1)^4|\bm{A}^{(8)}|$. Applying the hypothesis on
$\bm{A}^{(8)}$, one has
$|\bm{A}^{(8)}|=(-1)^{\frac{(m-1)(m-2)}{2}}\prod_{i=2}^mX_i\prod_{2\leq
i<j\leq m}(X_j-X_i)^4$ and then immediately gets
$|\bm{A}|=(-1)^{\frac{m(m-1)}{2}}\prod_{i=1}^mX_j\prod_{1\leq
i<j\leq m}(X_j-X_i)^4$.

From the above two cases, this theorem is thus proved.
\end{proof}

\subsection{Two Proofs of Theorem \ref{theorem:det-binom-powers}}

In this subsection, we give two inductive proofs of this theorem,
one uses induction on $m$ and another uses induction on $n$. The two
proofs are based on the same idea of reducing the matrix, though the
first proof is simpler in organization and understanding.

\subsubsection{The First Proof (Induction on $\bm{m}$)}

We first prove a lemma to simplify the first proof of Theorem
\ref{theorem:det-binom-powers}. This lemma is actually a special
case of the theorem under study when $m=1$ and $X_1=1$.

\begin{lemma}\label{lemma:det-binom}
When $1\leq m\leq n$, the determinant of the $m\times m$ matrix
$\bm{A}_{n,m}=\left[\binom{n+j-1}{i-1}\right]_{1\leq j\leq m \atop
1\leq i\leq m}$ is always equal to 1.
\end{lemma}
\begin{proof}
We use induction on $m$ to prove this lemma.

1) When $m=1$, $\bm{A}_{n,1}=[\binom{n}{0}]=[1]$. It is obvious that
$|\bm{A}_{n,1}|=1$.

2) Suppose this lemma is true for $m-1\geq 1$, let us prove the case
of $m\geq 2$. Write $\bm{A}_{n,m}$ as follows:
\begin{equation}
\left[\begin{matrix}%
1 & 1 & \cdots & 1\\
n & n+1 & \cdots & n+m-1\\
\binom{n}{2} & \binom{n+1}{2} & \cdots & \binom{n+m-1}{2}\\
\cdots & \cdots & \ddots & \cdots\\
\binom{n}{m-1} & \binom{n+1}{m-1} & \cdots & \binom{n+m-1}{m-1}
\end{matrix}\right].
\end{equation}
For $i=2\sim m$, subtract column $(i-1)$ column from column $i$, one
gets the following matrix:
\[
\left[\begin{matrix}%
1 & 0 & \cdots & 0\\
n & 1 & \cdots & 1\\
\binom{n}{2} & \binom{n+1}{2}-\binom{n}{2} & \cdots & \binom{n+m-1}{2}-\binom{n+m-2}{2}\\
\cdots & \cdots & \ddots & \cdots\\
\binom{n}{m-1} & \binom{n+1}{m-1}-\binom{n}{m-1} & \cdots &
\binom{n+m-1}{m-1}-\binom{n+m-2}{m-1}
\end{matrix}\right].
\]
From the property of binomial coefficients
\cite{Merris:Combinatorics2003},
$\binom{j}{i}-\binom{j-1}{i}=\binom{j-1}{i-1}$, so the above matrix
becomes
\[
\left[\begin{matrix}%
1 & 0 & \cdots & 0\\
n & 1 & \cdots & 1\\
\binom{n}{2} & \binom{n}{1} & \cdots & \binom{n+m-2}{1}\\
\cdots & \cdots & \ddots & \cdots\\
\binom{n}{m+1} & \binom{n}{m-2} & \cdots & \binom{n+m-2}{m-2}
\end{matrix}\right]=
\left[\begin{matrix}%
\begin{matrix}
1\\
n\\
\binom{n}{2}\\
\cdots\\
\binom{n}{m+1}
\end{matrix} &
\begin{matrix}
\begin{matrix}0 & 0 & \cdots & 0\end{matrix}\\
\\
\bm{A}_{n,m-1}
\\
\\
\\
\end{matrix}
\end{matrix}\right].
\]
Then, from the hypothesis, $|\bm{A}_{n,m}|=1\cdot
|\bm{A}_{n,m-1}|=1$. Thus this lemma is proved.
\end{proof}

\begin{proof}[The First Proof of Theorem \ref{theorem:det-binom-powers}]
In this proof, we use induction on $m$ to prove this theorem.

1) When $m=1$ and $n\geq l\geq 1$, $\bm{A}$ is simplified into an
$l\times l$ matrix as follows:
\[
\bm{A}=\left[\begin{matrix}%
\binom{n}{0} & \binom{n+1}{0}X_1 & \cdots & \binom{n+(l-1)}{0}X_1^{l-1}\\
\binom{n}{1} & \binom{n+1}{1}X_1 & \cdots & \binom{n+(l-1)}{1}X_1^{l-1}\\
\vdots & \vdots & \ddots & \vdots\\
\binom{n}{l-1} & \binom{n+1}{l-1}X_1 & \cdots &
\binom{n+(l-1)}{l-1}X_1^{l-1}
\end{matrix}\right]_{l\times l}.
\]
Factor the common terms in each column, the above matrix is reduced
to be
\[
\widetilde{\bm{A}}=\left[\begin{matrix}%
\binom{n}{0} & \binom{n+1}{0} & \cdots & \binom{n+(l-1)}{0}\\
\binom{n}{1} & \binom{n+1}{1} & \cdots & \binom{n+(l-1)}{1}\\
\vdots & \vdots & \ddots & \vdots\\
\binom{n}{l-1} & \binom{n+1}{l-1} & \cdots & \binom{n+(l-1)}{l-1}
\end{matrix}\right]_{l\times l}.
\]
From Lemma \ref{lemma:det-binom},
$\left|\widetilde{\bm{A}}\right|=1$, so
$|\bm{A}|=X_1^{1+\cdots+(l-1)}\left|\widetilde{\bm{A}}\right|=X_1^{\frac{l(l-1)}{2}}$,
which is equal to $\prod_{i=1}^1X_i^{\frac{l(l-1)}{2}}\prod_{1\leq
i<j\leq 1}(X_j-X_i)^{l^2}$ (the second term actually does not
exist).

2) Suppose this theorem is true for $m-1$ and $n\geq l\geq 1$, let
us prove the case of $m\geq 2$ and $n\geq l\geq 1$.

Before starting this part, we give a brief introduction to the basic
idea underlying the proof. The matrix $\bm{A}$ has a special feature
after the following elementary matrix operations: for $j=2\sim ml$,
subtracting column $j-1$ multiplied by $X_1$ from column $j$, row 1
of $\bm{A}$ becomes $[1\quad 0 \quad 0 \quad \cdots \quad 0]$. Then,
one can remove row 1 and column 1 from $\bm{A}$ and reduce $\bm{A}$
in some way. Repeat this process for $n$ rounds, $\bm{A}_1$ can be
completely removed from $\bm{A}$, which means that the value of $m$
decreases by one and the hypothesis can be applied to prove the
result of $m\geq 2$ and $n\geq l\geq 1$.

In the following, let us see how to reduce the matrix in the first
round of the process. Here, to achieve a clearer description of the
process, we use bracketed superscripts with increased digits to
denote the new matrices, each sub-matrices, and their elements after
different matrix operations (including reductions of the size). For
example, $\bm{A}^{(1)}$ denotes the matrix obtained after the above
subtractions, and
$\bm{A}_i^{(1)}=\left[a_{i,j,k}^{(1)}\right]_{1\leq j\leq ml \atop
1\leq k\leq l}$ denotes the $i$-th sub-matrix of $\bm{A}^{(1)}$.
Specially the original matrix is always written as $\bm{A}$ (without
any superscript) and its sub-matrix as $\bm{A}_i$.

For $j=2\sim ml$, multiplying column $j-1$ by $X_1$ and subtract it
from column $j$, the element of $\bm{A}_i^{(1)}$ at position $(j,1)$
becomes $a_{i,j,1}^{(1)}=X_i^{j-1}-X_i^{j-2}X_1=X_i^{j-2}(X_i-X_1)$
and the element at position $(j,k\geq 2)$ becomes:
\begin{eqnarray*}
a_{i,j,k}^{(1)} & = &
\binom{n+j-1}{k-1}X_i^{j-1}-\binom{n+j-2}{k-1}X_i^{j-2}X_1\\
& = &
\left(\binom{n+j-2}{k-1}+\binom{n+j-2}{k-2}\right)X_i^{j-1}-\binom{n+j-2}{k-1}X_i^{j-2}X_1\\
& = &
\binom{n+j-2}{k-1}X_i^{j-2}(X_i-X_1)+\binom{n+j-2}{k-2}X_i^{j-1}.
\end{eqnarray*}
When $i=1$, the above elements become: $a_{1,j,1}^{(1)}=0$ and
$a_{1,j,k}^{(1)}=\binom{n+j-2}{k-2}X_1^{j-1}$ ($k\geq 2$). So, row 1
of $\bm{A}_1$ become $[\begin{matrix}1 & 0 & 0 & \cdots &
0\end{matrix}]$. Then, $|\bm{A}|$ is equal to the determinant of the
following $(ml-1)\times(ml-1)$ matrix after removing row 1 and
column 1 of $\bm{A}^{(1)}$:
\[
\bm{A}^{(2)}=\left[\begin{matrix}%
\bm{A}_1^{(2)}\\
\bm{A}_2^{(2)}\\
\vdots\\
\bm{A}_m^{(2)}
\end{matrix}\right],
\]
where $\bm{A}_1^{(2)}=\left[\binom{n+j-1}{k-1}X_1^j\right]_{1\leq
j\leq mn-1 \atop 1\leq k\leq n-1}$ is an $(ml-1)\times(l-1)$ matrix,
and for $2\leq i\leq m$, $\bm{A}_1^{(2)}$ is an $(ml-1)\times l$
matrix as follows:
\[
\bm{A}_1^{(2)}=\left[a_{i,j,k}^{(2)}\right]_{1\leq j\leq ml-1 \atop
1\leq k\leq l}=\left[a_{i,j+1,k}^{(1)}\right]_{1\leq j\leq ml-1
\atop 1\leq k\leq l}=\left[\begin{matrix}%
X_i^{j-1}(X_i-X_1), & k=1\\
\binom{n+j-1}{k-1}X_i^{j-1}(X_i-X_1)+\binom{n+j-1}{k-2}X_i^j, &
k\geq 2
\end{matrix}\right]_{1\leq j\leq ml-1
\atop 1\leq k\leq l}.
\]
Then, let us reduce $\bm{A}^{(2)}$ to be of the same form as
$\bm{A}$. For $\bm{A}_1^{(2)}$, we simply factor out $X_1$ from each
row and get
$\bm{A}_1^{(3)}=\left[\binom{n+j-1}{k-1}X_1^{j-1}\right]_{1\leq
j\leq ml-1 \atop 1\leq k\leq l-1}$. Then, consider $\bm{A}_i^{(2)}$
for $i\geq 2$. Factoring out $(X_i-X_1)$ from row 1, one gets
$a_{i,j,1}^{(3)}=X_i^{j-1}$. Then, multiplying row 1 by $X_i$ and
subtracting it from row 2, one has
$a_{i,j,2}^{(3)}=\binom{n+j-1}{1}X_i^{j-1}(X_i-X_1)$. Then,
factoring out $(X_i-X_1)$ from row 2, one gets
$a_{i,j,2}^{(3)}=\binom{n+j-1}{1}X_i^{j-1}$. Repeat the above
procedure for other rows, one can finally get
$\bm{A}_i^{(3)}=\left[\binom{n+j-1}{k-1}X_i^{j-1}\right]_{1\leq
j\leq ml-1 \atop 1\leq k\leq l}$ and $(X_i-X_1)^n$ is factored out.
Combining the above results, we have
\[
|\bm{A}|=X_1^{l-1}\prod_{2\leq i\leq m}(X_i-X_1)^l|\bm{A}^{(3)}|.
\]
Note that the above equation becomes $|\bm{A}|=\prod_{2\leq i\leq
m}(X_i-X_1)^l|\bm{A}^{(3)}|$ when $l=1$. Observing the $m$
sub-matrices, one can see that each sub-matrix is of the same form
as the original one in $\bm{A}$, except that row $l$ and column $ml$
are removed from $\bm{A}_1$ and column $ml$ is removed from
$\bm{A}_i$ ($i\geq 2$).

Next, repeat the above process on $\bm{A}^{(3)}$, we can finally get
the following $(ml-2)\times(ml-2)$ matrix:
\[
\bm{A}^{(4)}=\left[\begin{matrix}%
\bm{A}_1^{(4)}\\
\bm{A}_2^{(4)}\\
\vdots\\
\bm{A}_m^{(4)}
\end{matrix}\right],
\]
where for
$\bm{A}_1^{(4)}=\left[\binom{n+j-1}{k-1}X_i^{j-1}\right]_{1\leq
j\leq ml-2 \atop 1\leq k\leq l-2}$, and for $2\leq i\leq m$,
$\bm{A}_i^{(4)}=\left[\binom{n+j-1}{k-1}X_i^{j-1}\right]_{1\leq
j\leq ml-2 \atop 1\leq k\leq l}$. In addition, we also have
\[
|\bm{A}^{(4)}|=X_1^{l-2}\prod_{2\leq i\leq
m}(X_i-X_1)^l|\bm{A}^{(3)}|,
\]
where note that $\bm{A}_1^{(3)}$ has only $l-1$ rows (one less than
$\bm{A}_1^{(1)}$).

Repeat the above procedure for $j=3\sim l$ rounds again, one can get
\[
|\bm{A}^{\langle j\rangle}|=X_1^{l-j}\prod_{2\leq i\leq
m}(X_i-X_1)^l|\bm{A}^{\langle j-1\rangle}|,
\]
where $\bm{A}^{\langle j\rangle}$ denotes the reduced matrix of size
$(ml-j)\times(ml-j)$ obtained after the $j$-th round of the above
process finishes, specially, $\bm{A}^{\langle
1\rangle}=\bm{A}^{(3)}$ and $\bm{A}^{\langle
2\rangle}=\bm{A}^{(4)}$.

After total $l$ rounds of the above process, one finally gets an
$(ml-l)\times(ml-l)$ matrix
\[
\bm{A}^{\langle n\rangle}=\left[\begin{matrix}%
\bm{A}_2^{\langle n\rangle}\\
\bm{A}_3^{\langle n\rangle}\\
\vdots\\
\bm{A}_m^{\langle n\rangle}
\end{matrix}\right],
\]
in which the first sub-matrix $\bm{A}_1$ is completely removed and
all other sub-matrices are untouched. Apparently, now
$\bm{A}^{\langle l\rangle}$ is a matrix of the same kind with
parameter $m-1$ and $l$.

Combining the relation between $|\bm{A}|$ and $\bm{A}^{\langle
1\rangle}$, and the relationships between $|\bm{A}^{\langle
j\rangle}|$ and $|\bm{A}^{\langle j-1\rangle}|$ ($2\leq j\leq l$),
one has
\begin{eqnarray}
|\bm{A}| & = & X_1^{(l-1)+\cdots+1}\prod_{2\leq i\leq
m}(X_i-X_1)^{l\cdot l}
\left|\bm{A}^{\langle l\rangle}\right|\nonumber\\
& = & X_1^{\frac{l(l-1)}{2}}\prod_{2\leq i\leq
m}(X_i-X_1)^{l^2}\left|\bm{A}^{\langle
l\rangle}\right|.\label{equation:inductive-m}
\end{eqnarray}

Then, applying the hypothesis for $\bm{A}^{\langle l\rangle}$, we
finally have
\begin{eqnarray*}
|\bm{A}| & = & \left(X_1^{\frac{l(l-1)}{2}}\prod_{2\leq i\leq
m}(X_i-X_1)^{l^2}\right)\cdot\left(\prod_{i=2}^mX_i^{\frac{l(l-1)}{2}}\prod_{2\leq
i<j\leq m}(X_j-X_i)^{l^2}\right)\\
& = & \prod_{i=1}^mX_i^{\frac{l(l-1)}{2}}\prod_{1\leq i<j\leq
m}(X_j-X_i)^{l^2}.
\end{eqnarray*}
This proves the case of $m\geq 2$ and $n\geq l\geq 1$.

From the above two cases, this theorem is thus proved.
\end{proof}

\subsubsection{The Second Proof (Induction on $\bm{l}$)}

\begin{proof}[The Second Proof of Theorem \ref{theorem:det-binom-powers}]
In this proof, we use induction on $l$ to prove this theorem.

1) When $l=1$, $m\geq 1$ and $n\geq l$, $\bm{A}$ is simplified into
an $m\times m$ matrix as follows:
\[
\bm{A}=\left[\begin{matrix}%
1 & X_1 & \cdots & X_1\\
1 & X_2 & \cdots & X_2^{m-1}\\
\vdots & \vdots & \ddots & \vdots\\
1 & X_m & \cdots & X_m^{m-1}
\end{matrix}\right]_{m\times m}.
\]
This is a Vandermonde matrix, so $|\bm{A}|=\prod_{1\leq i<j\leq
m}(X_j-X_i)$ for $m\geq 1$ \cite[\S4.4]{Zhang:MatrixTheory1999},
which is equal to
\[
\prod_{i=1}^mX_i^0\prod_{1\leq i<j\leq
m}(X_j-X_i)=\prod_{i=1}^mX_i^{\frac{1(1-1)}{2}}\prod_{1\leq i<j\leq
m}(X_j-X_i)^{1^2}.
\]

2) Suppose this theorem is true for $l-1$, $n\geq l$ and $m\geq 1$,
let us prove the case of $l\geq 2$, $n\geq l$ and $m\geq 1$.

Before starting this part, we give a brief introduction to the basic
idea underlying the proof. The matrix $\bm{A}$ has a special feature
after the following elementary matrix operations: for $i=1\sim m$
and $j=2\sim ml$, subtracting column $j-1$ multiplied by $X_i$ from
column $j$, row 1 of $\bm{A}_i$ becomes $[1\quad 0 \quad 0 \quad
\cdots \quad 0]$. After removing row 1 of each sub-matrix and column
1 of $\bm{A}$, the whole matrix is reduced to be of size
$m(l-1)\times m(l-1)$ and each sub-matrix is reduced to be of size
$m(l-1)\times(l-1)$. More importantly, after a series of matrix
operations, the matrix can be finally reduced to be a matrix of the
same form as the original one (with only different size). As a
result, we can then use the hypothesis on the case of $l-1$ and
$m,n$ to prove the result on $l$ and $m,n$.

As the first step, for $j=2\sim ml$, multiplying column $j-1$ by
$X_1$ and subtract it from column $j$, let us see how the matrix can
be reduced. In the following proof, to achieve a clearer description
of the process, we use bracketed superscripts with increased digits
to denote the new matrices, each sub-matrices, and their elements
after different matrix operations (including reductions of the
size). For example, $\bm{A}^{(1)}$ denotes the matrix obtained after
the above subtractions, and
$\bm{A}_i^{(1)}=\left[a_{i,j,k}^{(1)}\right]_{1\leq j\leq ml \atop
1\leq k\leq l}$ denotes the $i$-th sub-matrix of $\bm{A}^{(1)}$.
Specially the original matrix is always written as $\bm{A}$ (without
any superscript) and its sub-matrix as $\bm{A}_i$.

After the above subtraction transformations, the element of
$\bm{A}_i^{(1)}$ at position $(j,1)$ becomes
$a_{i,j,1}^{(1)}=X_i^{j-1}-X_i^{j-2}X_1=X_i^{j-2}(X_i-X_1)$ and the
element at position $(j,k\geq 2)$ becomes:
\begin{eqnarray*}
a_{i,j,k}^{(1)} & = &
\binom{n+j-1}{k-1}X_i^{j-1}-\binom{n+j-2}{k-1}X_i^{j-2}X_1\\
& = &
\left(\binom{n+j-2}{k-1}+\binom{n+j-2}{k-2}\right)X_i^{j-1}-\binom{n+j-2}{k-1}X_i^{j-2}X_1\\
& = &
\binom{n+j-2}{k-1}X_i^{j-2}(X_i-X_1)+\binom{n+j-2}{k-2}X_i^{j-1}.
\end{eqnarray*}
When $i=1$, the above elements become: $a_{1,j,1}^{(1)}=0$ and
$a_{1,j,k}^{(1)}=\binom{n+j-2}{k-2}X_1^{j-1}$ ($k\geq 2$). So, row 1
of $\bm{A}_1$ become $[\begin{matrix}1 & 0 & 0 & \cdots &
0\end{matrix}]$. Then, $|\bm{A}|$ is equal to the determinant of the
following $(ml-1)\times(ml-1)$ matrix after removing row 1 and
column 1 of $\bm{A}^{(1)}$:
\[
\bm{A}^{(2)}=\left[\begin{matrix}%
\bm{A}_1^{(2)}\\
\bm{A}_2^{(2)}\\
\vdots\\
\bm{A}_m^{(2)}
\end{matrix}\right],
\]
where $\bm{A}_1^{(2)}=\left[\binom{n+j-1}{k-1}X_1^j\right]_{1\leq
j\leq ml-1 \atop 1\leq k\leq l-1}$ is an $(ml-1)\times(l-1)$ matrix,
and for $2\leq i\leq m$, $\bm{A}_1^{(2)}$ is an $(ml-1)\times l$
matrix as follows:
\[
\bm{A}_1^{(2)}=\left[a_{i,j,k}^{(2)}\right]_{1\leq j\leq ml-1 \atop
1\leq k\leq l}=\left[a_{i,j+1,k}^{(1)}\right]_{1\leq j\leq ml-1
\atop 1\leq k\leq l}=\left[\begin{matrix}%
X_i^{j-1}(X_i-X_1), & k=1\\
\binom{n+j-1}{k-1}X_i^{j-1}(X_i-X_1)+\binom{n+j-1}{k-2}X_i^j, &
k\geq 2
\end{matrix}\right]_{1\leq j\leq ml-1
\atop 1\leq k\leq l}.
\]
Then, let us reduce $\bm{A}^{(2)}$ to be of the same form as
$\bm{A}$. For $\bm{A}_1^{(2)}$, we simply factor out $X_1$ from each
row and get
$\bm{A}_1^{(3)}=\left[\binom{n+j-1}{k-1}X_1^{j-1}\right]_{1\leq
j\leq ml-1 \atop 1\leq k\leq l-1}$. Then, consider $\bm{A}_i^{(2)}$
for $i\geq 2$. Factoring out $(X_i-X_1)$ from row 1, one gets
$a_{i,j,1}^{(3)}=X_i^{j-1}$. Then, multiplying row 1 by $X_i$ and
subtracting it from row 2, one has
$a_{i,j,2}^{(3)}=\binom{n+j-1}{1}X_i^{j-1}(X_i-X_1)$. Then,
factoring out $(X_i-X_1)$ from row 2, one gets
$a_{i,j,2}^{(3)}=\binom{n+j-1}{1}X_i^{j-1}$. Repeat the above
procedure for other rows, one can finally get
$\bm{A}_i^{(3)}=\left[\binom{n+j-1}{k-1}X_i^{j-1}\right]_{1\leq
j\leq ml-1 \atop 1\leq k\leq l}$ and $(X_i-X_1)^l$ is factored out.
Combining the above results, we have
\[
|\bm{A}|=X_1^{l-1}\prod_{2\leq i\leq m}(X_i-X_1)^l|\bm{A}^{(3)}|.
\]
Note that the above equation becomes
$|\bm{A}|=X_1^{l-1}|\bm{A}^{(3)}|$ when $m=1$. Observing the $m$
sub-matrices, one can see that each sub-matrix is of the same form
as the original one in $\bm{A}$, except that row $l$ and column $ml$
are removed from $\bm{A}_1$ and column $ml$ is removed from
$\bm{A}_i$ ($i\geq 2$).

Next, repeat the above process on $\bm{A}^{(3)}$ after replacing
$X_1$ by $X_2$. Due to the similarity of the whole process, we omit
the details and finally get the following $(ml-2)\times(ml-2)$
matrix:
\[
\bm{A}^{(4)}=\left[\begin{matrix}%
\bm{A}_1^{(4)}\\
\bm{A}_2^{(4)}\\
\vdots\\
\bm{A}_m^{(4)}
\end{matrix}\right],
\]
where for $1\leq i\leq 2$,
$\bm{A}_i^{(4)}=\left[\binom{n+j-1}{k-1}X_i^{j-1}\right]_{1\leq
j\leq ml-2 \atop 1\leq k\leq l-1}$, and for $3\leq i\leq m$,
$\bm{A}_i^{(4)}=\left[\binom{n+j-1}{k-1}X_i^{j-1}\right]_{1\leq
j\leq ml-2 \atop 1\leq k\leq l}$. In addition, we have
$|\bm{A}^{(4)}|=(-1)^{l-1}X_2^{l-1}(X_1-X_2)^{l-1}\prod_{3\leq i\leq
m}(X_i-X_2)^l|\bm{A}^{(3)}|$, where $(-1)^{l-1}$ is induced by the
fact that $a_{2,1,1}^{(3)}$ is at the position of $(l,1)$ in the
full matrix $\bm{A}^{(3)}$ (note that $\bm{A}_1^{(3)}$ has only
$l-1$ rows).

Repeat the above procedure for $j=3\sim m$, one can get
\[
|\bm{A}^{\langle j\rangle}|=(-1)^{(j-1)(l-1)}X_j^{l-1}\prod_{1\leq
i\leq j-1}(X_i-X_j)^{l-1}\prod_{j+1\leq i\leq
m}(X_i-X_j)^l|\bm{A}^{\langle j-1\rangle}|,
\]
where $\bm{A}^{\langle j\rangle}$ denotes the reduced matrix of size
$(ml-j)\times(ml-j)$ obtained after $j$ rounds of the above process,
specially, $\bm{A}^{\langle 1\rangle}=\bm{A}^{(3)}$ and
$\bm{A}^{\langle 2\rangle}=\bm{A}^{(4)}$. After total $m$ rounds of
the above process, one finally gets an $(ml-m)\times(ml-m)$ matrix
$\bm{A}^{\langle m\rangle}$, in which each sub-matrix is an
$(ml-m)\times(l-1)$ matrix defined by $\bm{A}_i^{\langle
m\rangle}=\left[\binom{n+j-1}{k-1}X_i^{j-1}\right]_{1\leq j\leq ml-m
\atop 1\leq k\leq l-1}$ ($i\geq 1$). Apparently, $\bm{A}^{\langle
m\rangle}$ is a matrix of the same kind with parameter $l-1$ and
$m,n$.

Combining the relation between $|\bm{A}|$ and $\bm{A}^{\langle
1\rangle}$, and the relationships between $|\bm{A}^{\langle
j\rangle}|$ and $|\bm{A}^{\langle j-1\rangle}|$ ($2\leq j\leq m$),
one has
\begin{eqnarray*}
|\bm{A}| & = &
(-1)^{(l-1)+\cdots+(m-1)(l-1)}\prod_{i=1}^mX_i^{l-1}\prod_{1\leq
i<j\leq m}(X_i-X_j)^{l-1}\prod_{1\leq
i<j\leq m}(X_j-X_i)^l\left|\bm{A}^{\langle m\rangle}\right|\\
& = & (-1)^{(l-1)+\cdots+(m-1)(nl-1)}\prod_{i=1}^mX_i^{l-1}\\
& & {}\cdot\left((-1)^{(l-1)+\cdots+(m-1)(l-1)}\prod_{1\leq i<j\leq
m}(X_j-X_i)^{l-1}\right)\prod_{1\leq i<j\leq
m}(X_j-X_i)^l\left|\bm{A}^{\langle m\rangle}\right|\\
& = & \prod_{i=1}^mX_i^{l-1}\prod_{1\leq i<j\leq
m}(X_j-X_i)^{2l-1}\left|\bm{A}^{\langle m\rangle}\right|.
\end{eqnarray*}

Then, applying the hypothesis on $\bm{A}^{\langle m\rangle}$, we
finally have
\begin{eqnarray*}
|\bm{A}| & = & \left(\prod_{i=1}^mX_i^{l-1}\prod_{1\leq i<j\leq
m}(X_j-X_i)^{2l-1}\right)\cdot\left(\prod_{i=1}^mX_i^{\frac{(l-1)(l-2)}{2}}\prod_{1\leq
i<j\leq m}(X_j-X_i)^{(l-1)^2}\right)\\
& = & \prod_{i=1}^mX_i^{\frac{l(l-1)}{2}}\prod_{1\leq i<j\leq
m}(X_j-X_i)^{l^2}.
\end{eqnarray*}
This proves the case of $l\geq 2$, $n\geq l$ and $m\geq 1$.

From the above two cases, this theorem is thus proved.
\end{proof}

\section*{Acknowledgments}

The author thank Dr. Arthur L.B. Yang with the Department of
Mathematics, Nankai University, for his Maple program of guessing
the formula of the determinant in Theorem
\ref{theorem:det-binom-powers}. The following friends and experts
are also appreciated for their help and comments on proving Theorem
\ref{theorem:det-binom-powers}: Prof. Guanrong Chen with the
Department of Electronic Engineering, City University of Hong Kong;
Dr. Jinhu L\"{u} and Dr. Hongbin Ma with the Institute of Systems
Science, Chinese Academy of Sciences; Dr. Zhong Li with the
Department of Computer Engineering, FernUniversit\"{a}t in Hagen,
Germany; Prof. Chi-Kwong Li with the Department of Mathematics, The
College of William \& Mary, USA; Dr. Jian Sun with the Microsoft
Research Asia and Dr. Tianshu Wang with the IBM China Research
Laboratory.

Specially, the author would like to thank Prof. Christian
Krattenthaler with the Camille Jordan Institute Universit\'{e}
Claude Bernard Lyon 1 (France) for pointing out that the two
determinants are actually special cases of Theorems 20 and 23 in his
work \cite{ADC1999}.

\bibliographystyle{unsrt}
\bibliography{determinant}

\end{document}